\documentclass[reqno]{amsart}
\usepackage{graphicx}

\usepackage{epstopdf}
\usepackage[justification=centering]{caption}
\usepackage{lineno,hyperref}
\usepackage{amsfonts,amsmath,latexsym,amssymb,amsthm}
\usepackage{geometry}
\newcommand{\hs}{{\mathcal H\,}}
\newcommand{\ws}{{\mathcal W\,}}
\newcommand{\vs}{{\mathcal V\,}}
\newcommand{\ls}{{\mathcal L\,}}
\modulolinenumbers[5]
\usepackage{makecell,multirow,diagbox}
\hypersetup{
	colorlinks=true,
	linkcolor=blue
}
\hypersetup{
	citecolor=blue
}
\newtheorem{definition}{Definition}
\newtheorem{theorem}{Theorem}
\newtheorem{lemma}{Lemma}
\usepackage[justification=centering]{caption}
\newtheorem{corollary}{Corollary}
\newtheorem{proposition}{Proposition}

\newtheorem{remark}{Remark}

\numberwithin{equation}{section}
\allowdisplaybreaks[4]

\begin{document}

\title{Fusion frames for operators and atomic systems}
\author{Yuxiang Xu}
\address{School of Mathematical Sciences, University of Electronic Science and Technology of China, 611731, P. R. China}
\email{scxuyuxiang@126.com}
\author{Dongwei Li}
\address{School of Mathematics, HeFei University of Technology, 230009, P. R. China}

\email{ldwmath@163.com}

\author{Jinsong Leng}
\address{School of Mathematical Sciences, University of Electronic Science and Technology of China, 611731, P. R. China}
\email{jinsongleng@126.com}




\begin{abstract}
Recently, fusion frames and frames for operators were considered as generalizations of frames in Hilbert spaces. In this paper, we generalize some of the known results in frame theory to fusion frames related to a linear bounded operator $K$ which we call $K$-fusion frames. We obtain new $K$-fusion frames by considering $K$-fusion frames with a class of bounded linear operators and construct new $K$-fusion frames from given ones. We also study the stability of $K$-fusion frames under small perturbations. We further give some characterizations of  atomic systems with subspace sequences.
\end{abstract}

\maketitle

\section{Introduction}
Frames were first introduced by Duffin and Schaeffer \cite{duffin1952class} to study some problems in nonharmonic Fourier series, reintroduced in 1986 by Daubechies, Grossmann, and Meyer \cite{daubechies1986painless}, and popularized from then on. Nice properties of frames make them very useful in characterization of function spaces and other fields of applications such as signal processing \cite{li2018frame}, sampling \cite{eldar2003sampling}, coding and communications \cite{han2018recovery}, filter bank theory \cite{bolcskei1998frame} and system modeling \cite{ward1998construction}.
Let $I$ be a countable index set and $\hs$ denote Hilbert space. A sequence $\{f_i\}_{i\in I}$ is called a frame for $\hs$ if there exist constants $0< A_1\le B_1< \infty$  such that
$$A_1\|f\|^2\le\sum_{i\in I}|\left\langle f,f_i\right\rangle |^2\le B_1\|f\|^2,~~~~{\rm{for ~all}}  ~~f\in\hs.$$

The fusion frame (frames of subspaces) which was introduced by Casazza and Kutyniok \cite{casazza2004frames}  is a natural generalization of frame theory and related to the construction of global frames from local frames in Hilbert spaces. Due to this property, fusion frames are special suiting for applications such as distributed processing, parallel processing of large frame systems \cite{casazza2008fusion}, optimal transmission by packet encoding \cite{bodmann2007optimal}.

Throughout this paper, $\hs$, $\hs_1$ and $\hs_2$ are separable Hilbert spaces. We denote by $\ls(\hs_1,\hs_2)$ the space of all bounded linear operators between $\hs_1$ and $\hs_2$. If $\hs_1=\hs_2=\hs$, then $\ls(\hs_1,\hs_2)$ is denoted by $\ls(\hs)$. For $T\in\ls(\hs)$, the range and the kernel of $T$ are denoted by $R(T)$ and $N(T)$, respectively, the pseudo-inverse of $T$ is denoted by $T^{\dagger}$.
\begin{definition}
	Let $\{W_i\}_{i\in I}$ be a sequence of closed subspaces in $\hs$, and let $\{w_i\}_{i\in I}$ be a family of weights, i.e., $w_i>0$ for all $i\in I$. We say that $\ws=\{(W_i,w_i)\}_{i\in I}$ is a fusion frame for $\hs$, if there exist constants $0<A_2\le B_2<\infty$ such that
	$$A_2\|f\|^2\le\sum_{i\in I}w_i^2\|\pi_{W_i}(f)\|^2\le B_2\|f\|^2,$$
	for all $f\in\hs$, where $\pi_{W_i}$ is the orthogonal projection onto the subspace $W_i$.
\end{definition}
We call $A_2$ and $B_2$ the fusion frame bounds. The family $\ws$ is called an $A_2$-tight fusion frame if $A_2=B_2$, it is a Parseval fusion frame if $A_2=B_2=1$. If $\ws=\{(W_i,w_i)\}_{i\in I}$ possesses an upper fusion frame bound, but not necessarily a lower, we say that $\ws$ is a Bessel fusion sequence with Bessel fusion bound $B_2$. Moreover, let $\{f_{ij}\}_{j\in J_i}$ be a frame for $W_i$ for each $i\in I$, then we call $\{(W_i,w_i,\{f_{ij}\}_{j\in J_i})\}_{i\in I}$ a fusion frame system for $\hs$.

For each Bessel fusion sequence $\ws=\{(W_i,w_i)\}_{i\in I}$ of $\hs$, we define the representation space associated with $\ws$ by
$$\bigg(\sum_{i\in I}\oplus W_i\bigg)_{\ell^2}=\bigg\{\{a_i\}_{i\in I}:a_i\in W_i,~\sum_{i\in I}\|a_i\|^2<\infty\bigg\}$$
with inner product given by

$$\langle \{a_i\}_{i\in I},\{c_i\}_{i\in I}\rangle =\sum_{i\in I}\langle a_i,c_i\rangle.$$

The following theorem provides a link between local and global properties.
\begin{theorem}\label{the1}
	{	\rm\cite{casazza2004frames}}
	For each $i\in I$, let $w_i>0$, let $W_i$ be a closed subspace of $\hs$, and let $\{f_{ij}\}_{j\in J_i}$ be a frame for $W_i$ with bounds $A_i$ and $B_i$. Suppose that $0<C=\inf_{i\in I} C_i\le \sup_{i\in I}D_i=D<\infty$. Then the following conditions are equivalent:
	\begin{enumerate}
		\item [\rm(i)] $\{(W_i,w_i)\}_{i\in I}$ is a fusion frame for $\hs$.
		\item [\rm(ii)]$\{w_if_{ij}\}_{j\in J_i,i\in I}$ is a frame for $\hs$.
	\end{enumerate}
	In particular, if $\{(W_i,w_i,\{f_{ij}\}_{j\in J_i})\}_{i\in I}$ is a fusion frame system for $\hs$ with fusion frame bounds $A$ and $B$, then $\{w_if_{ij}\}_{j\in J_i,i\in I}$ is a frame for $\hs$ with frame bounds $AC$ and $BD$.
\end{theorem}

The concept of $K$-frames was introduced by G{\u{a}}vru{\c{t}}a \cite{guavructa2012frames} to study the atomic systems \cite{feichtinger2001atomic}
with respect to a bounded linear operator $K$ in a separable Hilbert space $\hs$. It is known that $K$-frames are
more general than traditional frames in the sense that the lower frame bound only holds for the
elements in the range of $K$.
Several
methods to construct K-frames and the stability of perturbations for the K-frames have been
discussed in \cite{li2018generalized,xiao2013some}.

\begin{definition}
	Let $K$ be a bounded linear operator from $\hs$ to $\hs$. A sequence $\{f_i\}_{i\in I}$ in $\hs$ is called a $K$-frame for $\hs$ if there exist constants $0<A_3\le B_3<\infty$ such that
	$$A_3\|K^*f\|^2\le\sum_{i\in I}|\left\langle f,f_i\right\rangle |^2\le B_3\|f\|^2~~~~{\rm{for ~all}}  ~~f\in\hs.$$
\end{definition}

We call $A_3$, $B_3$ the lower frame bound and the upper frame bound for $K$-frame, respectively.

\begin{remark}
If $K=I_{\hs}$, then $K$-frames are just the ordinary frames, where $I_{\hs}$ is the identity operator on $\hs$.
\end{remark} 

Suppose that $\{f_i\}_{i\in I}$ is a $K$-frame for $\hs$. Obviously it is a Bessel sequence, then there exist three associated operators which are given as follows.

The synthesis operator is defined by
$$T:l^2(I)\rightarrow\hs,~~~T(c_i)=\sum_{i\in I}c_if_i.$$
Its adjoint operator
$$T^*:\hs\rightarrow l^2(I),~~~T^*f=\{\left\langle f,f_i\right\rangle \}_{i\in I}$$
is called the analysis operator. Then the frame operator is given by
$$S:\hs\rightarrow\hs,~~~Sf=\sum_{i\in I}\left\langle f,f_i\right\rangle f_i.$$
Note that, the frame operator of $K$-frame is not  invertible on $\hs$ in general and $S\ge AKK^*$.

Christensen and Heil \cite{christensen1997perturbations} gave the concept of atomic decompositions in Banach spaces. G{\u{a}}vru{\c{t}}a \cite{guavructa2012frames} gave the definition of atomic system for a bounded linear operator, and  obtained some characterizations of atomic systems as follows.
\begin{definition}
	{	\rm\cite{guavructa2012frames}}
	Let $K\in\ls(\hs)$. We say that $\{f_i\}_{i\in I}$ is an atomic system for $K$ if the following statements hold
	\begin{enumerate}
		\item[(i)]  the series $\sum_{i\in I}c_if_i$ converges for all $c=\{c_i\}_{i\in I}\in l^2$.
		\item[(ii)]  there exists $C>0$ such that for every $f\in\hs$ there exists $\{a_i\}_{i\in I}\in l^2$ such that $\|\{a_i\}\|_{l^2}\le C\|f\|$ and $Kf=\sum_{i\in I}a_if_i$.
	\end{enumerate}
\end{definition}
\begin{remark}\label{rem2}
The condition (i) in Definition 3 actually says that $\{f_i\}_{i\in I}$ is a Bessel sequence $($see Corollary $3.2.4$ of \cite{christensen2016introduction}$)$.
\end{remark} 
\begin{theorem}
	{	\rm\cite{guavructa2012frames}}
	Let $\{f_i\}_{i\in I}\subset\hs$ and $K$ be a bounded linear operator.  $\{f_i\}_{i\in I}$ is an atomic system for $K$ if and only if $\{f_i\}_{i\in I}$ is a $K$-frame for $\hs$.
\end{theorem}
\begin{lemma}
	{	\rm{\cite{christensen2016introduction}}}
	Let $\hs$ be a Hilbert space, and suppose that $T\in \ls({\hs})$ has a closed range. Then there exists an operator $T^{\dagger}\in \ls({\hs})$ for which
	$$N(T^{\dagger})=R(T)^{\perp},~~R(T^{\dagger})=N(T)^{\perp},~~TT^{\dagger}y=y,~~y\in R(T).$$
	We call the operator $T^{\dagger}$ the pseudo-inverse of $T$. This operator is uniquely determined by these properties.
\end{lemma}
In fact, if $T$ is invertible, then we have $T^{-1}=T^{\dagger}$.
\begin{theorem}\label{thm3}
	{	\rm\cite{douglas1966majorization}}
	(Douglas' majorization theorem ). Let $U\in\ls(\hs_1,\hs)$, $V\in\ls(\hs_2,\hs)$ be two bounded operators. Then the following are equivalent:
	\begin{enumerate}
		\item [(i)] R(U)$\subset$R(V);
		\item [(ii)] $UU^{*}\le \lambda^2VV^{*}$ for some $\lambda> 0$ (i.e., $V^{*}$ majorizes $U^{*}$);
		\item [(iii)]  $U=VT$ for some $T\in\ls(\hs_1,\hs_2)$.
	\end{enumerate}
\end{theorem}

\begin{lemma}\label{lem2}
	{	\rm\cite{guavructa2007duality}}
	Let $T\in\ls(\hs)$ and $W$ be a closed subspace of $\hs$. Then we have
	$$\pi_{W}T^{*}=\pi_{W}T^*\pi_{\overline{TW}}.$$
\end{lemma}
\section{Fusion frames for operators}
In this section, we generalize some of the known results in $K$-frame theory to $K$-fusion frames.
\begin{definition}\label{def4}
	Let $\{W_i\}_{i\in I}$ be a sequence of closed subspaces in $\hs$, and let $\{w_i\}_{i\in I}$ be a family of weights, i.e., $w_i>0$ for all $i\in I$. Let $K\in\ls(\hs)$, we say that $\ws=\{(W_i,w_i)\}_{i\in I}$ is a $K$-fusion frame for $\hs$, if there exist constants $0<A\le B<\infty$ such that
	\begin{equation}\label{2.1}
	A\|K^*f\|^2\le\sum_{i\in I}w_i^2\|\pi_{W_i}(f)\|^2\le B\|f\|^2,
	\end{equation}
	for all $f\in\hs$, where $\pi_{W_i}$ is the orthogonal projection onto the subspace $W_i$.
\end{definition}

$A$, $B$ are called the $K$-fusion frame bounds, respectively. If $w_i=w_j=w$ for all $i,j\in I$, we call $\ws$ a $w$-uniform $K$-fusion frame for $\hs$. Moreover, we say that $\ws=\{(W_i,w_i)\}_{i\in I}$ is an orthonormal $K$-fusion basis for $\hs$ if $\hs=\bigoplus_{i\in I}W_i$.

\begin{remark}
 Let $\ws=\{(W_i,w_i)\}_{i\in I}$ be a $K$-fusion frame for $\hs$. And let $W_i={\rm span}\{f_{ij}\}_{j\in J_i}$  for each $i\in I$. If $\{f_{ij}\}_{j\in J_i}$ is a frame for $W_i$, then we call
$\ws=\{(W_i,w_i,\{f_{ij}\}_{j\in J_i})\}_{i\in I}$ a $K$-fusion frame system for $\hs$. Moreover, $\{w_if_{ij}\}_{j\in J_i,i\in I}$ is a $K$-frame for $\hs$.
\end{remark}

Let $\ws$ be a Bessel fusion sequence for $\hs$. The synthesis operator is defined by
$$T_{\ws}:\bigg(\sum_{i\in I}\oplus W_i\bigg)_{\ell^2}\rightarrow\hs,~~T_{\ws}(\{a_i\})=\sum_{i\in I}w_ia_i,~~\forall a=\{a_i\}_{i\in I}\in \bigg(\sum_{i\in I}\oplus W_i\bigg)_{\ell^2}.$$
Similar to the fusion frame in Hilbert space, $T_{\ws}(\{a_i\})$ converges unconditionally (see Lemma 3.9 of \cite{casazza2004frames}).
The adjoint operator $T^{*}_{\ws}:\hs\rightarrow \bigg(\sum_{i\in I}\oplus W_i\bigg)_{\ell^2}$ given by $T^{*}_{\ws}(f)=\{w_i\pi_{W_i}(f)\}_{i\in I}$ is called  the analysis operator.

The $K$-fusion frame operator $S_{\ws}$ for $\ws$ is defined by
$$S_{\ws}:\hs\rightarrow\hs,~~S_{\ws}(f)=\sum_{i\in I}w_i^2\pi_{W_i}(f).$$

\begin{proposition}\label{prop1}
	Let $K\in\ls(\hs)$ be with closed range. Let $\ws=\{(W_i,w_i)\}_{i\in I}$ be a $K$-fusion frame for $\hs$ with frame operator $S_{\ws}$, then the following holds,
	\begin{enumerate}
		\item[(i)] $AKK^{*}\le S_{\ws}\le B\cdot I_{\hs}$;
		\item[(ii)] $A\|K^{\dagger}\|^{-2}\|f\|\le\|S_{\ws}f\|\le B\|f\|,~~\forall f\in R(K)$;
		\item[(iii)]  $B^{-1}\|f\|\le\|S^{-1}_{\ws}f\|\le A^{-1}\|K^{\dagger}\|^2\|f\|,~~\forall f\in S_{\ws}(R(K))$.
	\end{enumerate}
\end{proposition}

\begin{proof}

 (i).  Let $\ws=\{(W_i,w_i)\}_{i\in I}$ be a $K$-fusion frame for $\hs$ with frame operator $S_{\ws}$, then
$$A\|K^*f\|^2\le\sum_{i\in I}w_i^2\|\pi_{W_i}(f)\|^2=\left\langle S_{\ws}f,f\right\rangle \le B\|f\|^2,$$
that is,
\begin{equation}\label{2.2}
	\left\langle AKK^*f,f\right\rangle \le \left\langle S_{\ws}f,f\right\rangle\le \left\langle Bf,f\right\rangle ,~~\forall f\in\hs.
\end{equation}
(ii). Since $R(K)$ is closed, there exists a pseudo-inverse $K^{\dagger}$ of $K$ such that $KK^{\dagger}f=f$, $\forall f\in R(K)$, namely $KK^{\dagger}|_{R(K)}=I_{R(K)}$, so we have
$$I^*_{R(K)}=(K^{\dagger}|_{R(K)})^*K^*.$$
Hence for any $f\in R(K)$, we have
$$\|f\|=\|(K^{\dagger}|_{R(K)})^*K^*f\|\le \|K^{\dagger}\|\|K^*f\|,$$
that is, $$\|K^*f\|^{2}\ge \|K^{\dagger}\|^{-2}\|f\|^2.$$
Combined with (\ref{2.2}) we have
$$A\|K^{\dagger}\|^{-2}\|f\|^2\le A\|K^*f\|^2\le \left\langle S_{\ws}f,f\right\rangle,~~\forall f\in R(K).$$
So, from the Definition 4  we have
\begin{equation}\label{2.3}
A\|K^{\dagger}\|^{-2}\|f\|\le \| S_{\ws}f\|\le B\|f\|,~~\forall f\in R(K).
\end{equation}
(iii). In terms of (\ref{2.3}), we have
$$B^{-1}\|f\|\le\|S^{-1}_{\ws}f\|\le A^{-1}\|K^{\dagger}\|^2\|f\|,~~\forall f\in S_{\ws}(R(K)).$$
\end{proof}

Next, we use the frame operator of $K$-fusion frame to give another equivalent characterization of $K$-fusion frame. First, we need the following important result from operator theory:

\begin{proposition}\label{prop2}
	Let  $\ws=\{(W_i,w_i)\}_{i\in I}$ be a Bessel fusion sequence in $\hs$. Then $\ws$ is a $K$-fusion frame for $\hs$, if and only if there exists $A>0$ such that $AKK^*\le S_{\ws}$, where $S_{\ws}$ is the frame operator for $\ws$.
\end{proposition}
\begin{proof}
 By the Definition \ref{def4},  $\ws=\{(W_i,w_i)\}_{i\in I}$ is a $K$-fusion frame for $\hs$ with bounds $A,B$ and frame operator $S_{\ws}$, if and only if
$$A\|K^*f\|^2\le\sum_{i\in I}v_i\|\pi_{W_i}(f)\|^2=\left\langle S_{\ws}f,f\right\rangle \le B\|f\|^2,~~~\forall f\in\hs. $$
Therefore,
$$\left\langle AKK^*f,f\right\rangle \le \left\langle S_{\ws}f,f\right\rangle\le \left\langle Bf,f\right\rangle, ~~~\forall f\in\hs.$$
And so by Proposition \ref{prop1} the  conclusion holds.
\end{proof}
\begin{corollary}
	Let  $\ws=\{(W_i,w_i)\}_{i\in I}$ be a $K$-fusion frame for $\hs$. Let $T\in\ls(\hs)$ with $R(T)\subset R(K)$. Then $\ws$ is a $T$-fusion frame for $\hs$.
\end{corollary}
\begin{proof} This follows immediately from Proposition \ref{prop2}. 
\end{proof}

\begin{corollary}
	Let  $\ws=\{(W_i,w_i)\}_{i\in I}$ be a Bessel fusion sequence in $\hs$ with frame operator $S_{\ws}$. Then $\ws$ is a $K$-fusion frame for $\hs$  if and only if $K=S_{\ws}^{1/2}T$, for some $T\in\ls(\hs)$.
\end{corollary}
\begin{proof}
By Theorem 4, $\ws$ is a $K$-fusion frame if and only if there exists $A>0$ such that
$$AKK^*\le S_{\ws}=S^{1/2}_{\ws}{S^{1/2}_{\ws}}^*.$$
Therefore by Proposition \ref{prop2} the conclusion hold. 
\end{proof}
In the following we give several characterizations when a bounded linear operator $T$ is applied to a $K$-fusion frame.

\begin{theorem}
	Let $K\in\ls(\hs)$ be with a dense range. Let $\ws=\{(W_i,w_i)\}_{i\in I}$ be a $K$-fusion frame and $T\in\ls(\hs)$ have closed range. If $\{(TW_i,w_i)\}_{i\in I}$ is a $K$-fusion frame for $\hs$, then $T$ is surjective.
\end{theorem}
\begin{proof} Suppose $\{(TW_i,w_i)\}_{i\in I}$ is a $K$-fusion frame for $\hs$ with frame bounds $A$ and $B$. Suppose  $\{f_{ij}\}_{j\in J_i}$ is a frame for $W_i$ with bounds $C_i$ and $D_i$, where $0<C=\inf_{i\in I} C_i\le \sup_{i\in I}D_i=D<\infty$. From  Theorem \ref{the1}, $\{Tw_if_{ij}\}_{j\in J_i,i\in I}$ is a $K$-frame for $\hs$ with bounds $AC$ and $BD$. Then for any $f\in\hs$,
	\begin{equation}\label{2.4}
AC\|K^*f\|^2\le\sum_{i\in I}\sum_{j\in J_i}|\langle f,Tw_if_{ij}\rangle|^2\le BD\|f\|^2.
	\end{equation}
As $K$ is with a dense rang, $K^*$ is injective. Then from (\ref{2.4}), $T^*$ is injective since $N(T^*)\subset N(K^*)$. Moreover,  $R(T)=N(T^*)^{\perp}=\hs$. Thus $T$ is surjective.
\end{proof}
\begin{theorem}\label{thm5}
	Let $K\in\ls(\hs)$ and $\ws=\{(W_i,w_i)\}_{i\in I}$ be a $K$-fusion frame for $\hs$. Let $T\in\ls(\hs)$ has closed range with $TK=KT$ and $T^{\dagger}T(W_i)\subset W_i$. If $T$ is surjective, then $\{(TW_i,w_i)\}_{i\in I}$ is a $K$-fusion frame for $\hs$, where $T^{\dagger}$ is the pseudo-inverse of $T$.
\end{theorem}

\begin{proof}
Since $T$ has closed range, it has the  pseudo-inverse  $T^{\dagger}$ such that $TT^{\dagger}=I_{\hs}$. We now prove that $T(W_i)$ is again a closed subspace. Since $T^{\dagger}\overline{T(W_i)}\subset W_i$ and $T$ is surjective, we obtain $TT^{\dagger}\overline{T(W_i)}\subset T(W_i)$.
By Lemma 2.5.2 of \cite{christensen2016introduction}, we have $T^{\dagger}=T^*(TT^*)^{-1}$.
Hence
$$TT^*(TT^*)^{-1}\overline{T(W_i)}\subset T(W_i),$$
and then $\overline{T(W_i)}\subset T(W_i)$. Therefore, $T(W_i)$ is again a closed subspace.

From Lemma \ref{lem2} we have
$$\|\pi_{W_i}T^*f\|=\|\pi_{W_i}T^*\pi_{TW_i}f\|\le \|\pi_{W_i}\|\|T^*\|\|\pi_{TW_i}f\|\le\|T^*\|\|\pi_{TW_i}f\| ,$$
then
\begin{equation}\label{2.5}
\|\pi_{TW_i}f\|\ge \|{T}^*\|^{-1}\|\pi_{W_i}T^*f\|.
\end{equation}

Let $I_{\hs}=I_{\hs}^{*}={T^{\dagger}}^*T^*$, then for each $f\in R(T)=\hs$, $K^*f={T^{\dagger}}^*T^*K^*f$, so we have
$$\|K^*f\|=\|{T^{\dagger}}^*T^*K^*f\|\le \|{T^{\dagger}}^*\|\cdot\|T^*K^*f\|.$$
Since $TK=KT$, we have
\begin{equation}\label{2.6}
\|{T^{\dagger}}^*\|^{-1}\|K^*f\|\le \|T^*K^*f\|=\|K^*T^*f\|.
\end{equation}

Now for each $f\in \hs$, we obtain
\begin{eqnarray}
\sum_{i\in I}w_i^2\|\pi_{TW_i}f\|^2&\ge& \|{T}^*\|^{-2}\sum_{i\in I}w_i^2\|\pi_{W_i}T^*f\|^2  \nonumber\\
&\ge& \|{T}^*\|^{-2}A\|K^*T^*f\|^2\nonumber\\
&\ge&A \|{T}^*\|^{-2}\|{T^{\dagger}}^*\|^{-2} \|K^*f\|^2.   \nonumber
\end{eqnarray}
On the other hand, from Lemma \ref{lem2}, we obtain, with $T^{\dagger}$ instead of $T$:
$$\pi_{TW_i}{T^{\dagger}}^*=\pi_{TW_i}{T^{\dagger}}^{*}\pi_{W_i},$$
hence
\begin{equation}\label{2.7}
\|\pi_{TW_i}f\|=\|\pi_{TW_i}{T^{\dagger}}^{*}\pi_{W_i}T^*f\|\le \|\pi_{TW_i}\|\|{T^{\dagger}}^{*}\|\|\pi_{W_i}T^*f\|\le\|{T^{\dagger}}^{*}\|\|\pi_{W_i}T^*f\|.
\end{equation}

Since $\{(TW_i,w_i)\}_{i\in I}$ is a Bessel fusion sequence with bound $B$, from (\ref{2.7}), we have
\begin{eqnarray}
\sum_{i\in I}w_i^2\|\pi_{TW_i}f\|^2
&\le& \|{T^{\dagger}}^{*}\|^2\sum_{i\in I}w_i^2\|\pi_{W_i}T^*f\|^2  \nonumber\\
&\le& \|{T^{\dagger}}^{*}\|^2B\|T^*f\|^2\nonumber\\
&\le&B \|{T^{\dagger}}^*\|^{2}\|T\|^2 \|f\|^2.   \nonumber
\end{eqnarray}
Therefore, $\{(TW_i,w_i)\}_{i\in I}$ is a $K$-fusion frame for $\hs$.
\end{proof}

\begin{theorem}\label{thm6}
	Let $K\in\ls(\hs)$ be with a dense range. Let $\ws=\{(W_i,w_i)\}_{i\in I}$ be a $K$-fusion frame and let $T\in\ls(\hs)$ be with $TK = KT$. If $\{(TW_i,w_i)\}_{i\in I}$ and $\{(T^*W_i,w_i)\}_{i\in I}$ are $K$-fusion frames for $\hs$, then $T$ is invertible.
\end{theorem}

\begin{proof} Since $\{(TW_i,w_i)\}_{i\in I}$ is a $K$-fusion frame, from Theorem \ref{thm5}, we have $T$ is surjective.

Suppose $\{(T^*W_i,w_i)\}_{i\in I}$ is a $K$-fusion frame for $\hs$ with frame bounds $A$ and $B$. Then for any $f\in\hs$,
$$A\|K^*f\|^2\le\sum_{i\in I}w_i^2\|\pi_{T^*W_i}f\|^2\le B\|f\|^2.$$
Then $N(T)\subset N(K^*)$. As $K$ is with a dense range, $K^*$ is injective. Then we have $T$ is  injection. Therefore, $T$ is bijective.
By Bounded Inverse Theorem, $T$ is invertible.
\end{proof}

\begin{corollary}
	Let $\ws=\{(W_i,w_i)\}_{i\in I}$ be a $K$-fusion frame and let $T\in\ls(\hs)$ have closed range with $TK=KT$ and $T^*$ be isometry. Let $T^{\dagger}$ be the pseudo-inverse of $T$ and satisfy $T^{\dagger}T(W_i)\subset W_i$, then $\{(TW_i,w_i)\}_{i\in I}$ is a $K$-fusion frames for $R(T)$.
\end{corollary}
\begin{proof}
Since $T^*$ is isometry, for any $f\in\hs$, then we have
\begin{equation}\label{2.8}
\|K^*T^*f\|^2=\|T^*K^*f\|^2=\|K^*f\|^2.
\end{equation} 

Similar to the proof of Theorem \ref{thm6}, with (\ref{2.8}) instead of (\ref{2.6}), we can easy obtain the conclusion.
\end{proof}

The following result shows a class of operator in $\ls(\hs)$ associated with a given $K$-fusion frame.
\begin{theorem}
	Let $K_1,~K_2\in\ls(\hs)$ and $\ws=\{(W_i,w_i)\}_{i\in I}$ satisfy $(\ref{2.1})$ with $K_1$ and $K_2$. For any $\alpha,~\beta\in\mathbb{R}$ , then $\ws$ is a   $(\alpha K_1+\beta K_2)$-fusion frame $($and $(K_1K_2)$-fusion frame$)$ .
\end{theorem}
\begin{proof}
	Since $\ws=\{(W_i,w_i)\}_{i\in I}$ is an $K_n$-fusion frame for $n=1,~2$, there are positive constants $A_n,~B_n>0~(n=1,2)$ such that 
	\begin{eqnarray}\label{2.9}
	A_n\|K^*_nf\|^2\le \sum_{i\in I}w_i^2\|\pi_{W_i}(f)\|^2\le B_n\|f\|^2,~~~\forall f\in\hs.
	\end{eqnarray}
	Since 
	\begin{align*}
	\|K_1^*f\|^2&=\frac{1}{|\alpha|^2}\|\alpha K_1^*f\|^2=\frac{1}{|\alpha|^2}\|(\alpha K _1^*+\beta K^*_2)f-\beta K^*_2f\|^2\\
	&\ge \frac{1}{|\alpha|^2}\|(\alpha K _1^*+\beta K^*_2)f\|^2-\frac{1}{|\alpha|^2}\|\beta K^*_2f\|^2,
	\end{align*}
	we have 
	\begin{align*}
	\|(\alpha K _1^*+\beta K^*_2)f\|^2&\le |\alpha|^2\|K_1^*f\|^2+|\beta|^2\|K_2^*f\|^2\\
	&\le \frac{1}{2}(|\alpha|^2\|K_1^*f\|^2+|\beta|^2\|K_2^*f\|^2\\
	&\qquad+\frac{A_1}{A_2}|\beta|^2\|K_1^*f\|^2+\frac{A_2}{A_1}|\alpha|^2\|K_2^*f\|^2)\\
	&=\frac{A_2|\alpha|^2+A_1|\beta|^2}{2A_1A_2}(A_1\|K_1^*f\|^2+A_2\|K_2^*f\|^2).
	\end{align*}
	Hence
	$$\sum_{i\in I}w_i^2\|\pi_{W_i}(f)\|^2\ge \frac{1}{2}(A_1\|K_1^*f\|^2+A_2\|K_2^*f\|^2)\ge \frac{A_1A_2}{A_2|\alpha|^2+A_1|\beta|^2}\|(\alpha K _1^*+\beta K^*_2)f\|^2.$$
	And from inequalities \eqref{2.9}, we get 
	$$\sum_{i\in I}\sum_{i\in I}w_i^2\|\pi_{W_i}(f)\|^2\le \frac{B_1+B_2}{2}\|f\|^2,~~~\forall f\in\hs.$$
	Therefore, $\ws$ is a $(\alpha K_1+\beta K_2)$-fusion frame.
	
	Now for each $f\in\hs$, we have
	$$\|(K_1K_2)^*f\|^2=\|K^*_2K^*_1f\|^2\le \|K^*_2\|^2\|K_1^*f\|^2.$$
	From \eqref{2.1},
	$$\frac{A_1}{\|K^*_2\|^2}\|(K_1K_2)^*f\|^2\le A_1\|K^*_1f\|^2\le \sum_{i\in I}\|\Lambda_if\|^2\le B_1\|f\|^2,~~~\forall f\in\hs.$$
	Thus,  $\ws$ is a $(K_1K_2)$-fusion frame for $\hs$.
\end{proof}
Next, we construct new $K$-fusion frame with given one.
\begin{theorem}
	Let $K\in\ls(\hs)$ with closed range, $\ws=\{(W_i,w_i)\}_{i\in I}$ be a $K$-fusion frame for $\hs$ with bounds $A,~B$. Let $\vs=\{(V_i,v_i)\}_{i\in I}$ be a sequence with synthesis operator $T_{\vs}$. For any two positively confined sequences $a=\{a_i\}_{i\in I}$ and $b=\{b_i\}_{i\in I}$, if $\|T_{\vs}\|^2<\frac{A\inf_{i\in I}a_i^2}{2\|K^{\dagger}\|^2\sup_{i\in I}b_i^2}$, then $a\ws+b\vs$ is a $K$-fusion frame for $\hs$.
\end{theorem}
\begin{proof}
	For any $f\in\hs$, we have
	\begin{align*}
	&\sum_{i\in I}(a_i^2w_i^2\|\pi_{W_i}(f)\|^2+b_i^2v_i^2\|\pi_{V_i}(f)\|^2)\\
	&= \sum_{i\in I}a_i^2w_i^2\|\pi_{W_i}(f)\|^2+\sum_{i\in I}b_i^2v_i^2\|\pi_{V_i}(f)\|^2+2{\rm Re}\sum_{i\in I}\left\langle a_iw_i\pi_{W_i}(f),b_iv_i\pi_{V_i}(f)\right\rangle \\
	&\le 2(\sum_{i\in I}a_i^2w_i^2\|\pi_{W_i}(f)\|^2+\sum_{i\in I}b_i^2v_i^2\|\pi_{V_i}(f)\|^2\\
	&\le 2((\sup_{i\in I}a_i^2)\sum_{i\in I}w_i^2\|\pi_{W_i}(f)\|^2+(\sup_{i\in I}b_i^2)\sum_{i\in I}v_i^2\|\pi_{V_i}(f)\|^2)\\
	&\le 2((\sup_{i\in I}a_i^2)B\|f\|^2+(\sup_{i\in I}b_i^2)\|T_{\vs}^*f\|^2)\\
	&\le 2((\sup_{i\in I}a_i^2)B+(\sup_{i\in I}b_i^2)\|T_{\vs}\|^2)\|f\|^2.
	\end{align*}
	Since 
	\begin{align*}
	\sum_{i\in I}\|a_iw_i\pi_{W_i}(f)\|^2&=\sum_{i\in I}\|(a_iw_i\pi_{W_i}+b_iv_i\pi_{V_i})f-b_iv_i\pi_{V_i}f\|^2\\
	&\le 2(\sum_{i\in I}\|(a_iw_i\pi_{W_i}+b_iv_i\pi_{V_i})f\|^2+\sum_{i\in I}b_i^2v_i^2\|\pi_{V_i}(f)\|^2),
	\end{align*}
	we have
	\begin{align*}
2\sum_{i\in I}\|(a_iw_i\pi_{W_i}+b_iv_i\pi_{V_i})f\|^2&\ge\sum_{i\in I}\|a_iw_i\pi_{W_i}(f)\|^2-2\sum_{i\in I}b_i^2v_i^2\|\pi_{V_i}(f)\|^2\\
	&\ge(\inf_{i\in I}a_i^2)\sum_{i\in I}w_i^2\|\pi_{W_i}(f)\|^2-2(\sup_{i\in I}b_i^2)\|T_{\vs}^*f\|^2\\
	&\ge(\inf_{i\in I}a_i^2)\|K^*f\|^2-2(\sup_{i\in I}b_i^2)\|T_{\vs}^*(K^*)^{\dagger}f\|^2\\
	&\ge((\inf_{i\in I}a_i^2)-2(\sup_{i\in I}b_i^2)\|T_{\vs}\|^2\|K^{\dagger}\|^2)\|K^*f\|^2
	\end{align*}
	From $\|T_{\vs}\|^2<\frac{A\inf_{i\in I}a_i^2}{2\|K^{\dagger}\|^2\sup_{i\in I}b_i^2}$,
	we obtain that  $a\ws+b\vs$ is a $K$-fusion frame for $\hs$. 
\end{proof}
We give a sufficient condition under which a g-frame is stable under small perturbations in terms of positively confined sequence.
\begin{theorem}
	Let $\{(W_i,w_i)\}_{i\in I}$ be a $K$-fusion frame with bounds $A,~B$ for $\hs$ and let $\{a_i\}_{i\in I}$ and $\{b_i\}_{i\in I} \subset \mathbb{R}$ be two positively confined sequences. Suppose that   $V_i\subset \hs$ and there exist constants $0\le \lambda_1,~\lambda_2< 1/2$, $\mu>0 $ such that  
	$$\sum_{i\in I}w_i^2\|(a_i\pi_{W_i}-b_i\pi_{V_i})f\|^2\le \lambda_1\sum_{i\in I}a_i^2w_i^2\|\pi_{W_i}f\|^2+\lambda_2\sum_{i\in I}b_i^2w_i^2\|\pi_{V_i}f\|^2+\mu\|K^*f\|^2,~~f\in\hs.$$
	Then $\{(V_i,w_i)\}_{i\in I}$ is a $K$-fusion frame for $\hs$ with bounds
$$\frac{(1-2\lambda_1)(\inf_{i\in I}a_i)^2}{2(1+\lambda_2+\frac{\mu}{A})(\sup_{i\in I}b_i)^2}A~ {\rm and}~\frac{2(1+\lambda_1+\frac{\mu}{A})(\sup_{i\in I}a_i)^2}{(1-2\lambda_2)(\inf_{i\in I}b_i)^2}B.$$
\end{theorem}
\begin{proof}
	For any $f\in\hs$,
	\begin{align*}
	\sum_{i\in I}w_i^2\|b_i\pi_{V_i}f\|^2=&\sum_{i\in I}w_i^2\|(b_i\pi_{V_i}-a_i\pi_{W_i})f+a_i\pi_{W_i}f\|^2\\
	&\le 2(\sum_{i\in I}w_i^2\|(a_i\pi_{W_i}-b_i\pi_{V_i})f\|^2+\sum_{i\in I}w_i^2\|a_i\pi_{W_i}f\|^2)\\
	&\le  2(\lambda_1\sum_{i\in I}a_i^2w_i^2\|\pi_{W_i}f\|^2+\lambda_2\sum_{i\in I}b_i^2w_i^2\|\pi_{V_i}f\|^2+\mu\|K^*f\|^2+\sum_{i\in I}w_i^2\|a_i\pi_{W_i}f\|^2),
	\end{align*}
	then we have 
	\begin{align*}
	(1-2\lambda_2)(\inf_{i\in I}b_i)^2\sum_{i\in I}w_i^2\|\pi_{V_i}f\|^2\le 2(1+\lambda_1+\frac{\mu}{A})(\sup_{i\in I}a_i)^2\sum_{i\in I}w_i^2\|\pi_{W_i}f\|^2.
	\end{align*}
	Thus
	$$\sum_{i\in I}w_i^2\|\pi_{V_i}f\|^2\le \frac{2(1+\lambda_1+\frac{\mu}{A})(\sup_{i\in I}a_i)^2}{(1-2\lambda_2)(\inf_{i\in I}b_i)^2}\sum_{i\in I}\|\Lambda_if\|^2\le \frac{2(1+\lambda_1+\frac{\mu}{A})(\sup_{i\in I}a_i)^2}{(1-2\lambda_2)(\inf_{i\in I}b_i)^2}B\|f\|^2,~~\forall f\in\hs.$$
	On the other hand, since
	\begin{align*}
	\sum_{i\in I}w_i^2\|a_i\pi_{W_i}f\|^2=&\sum_{i\in I}w_i^2\|(a_i\pi_{W_i}-b_i\pi_{V_i})f+b_i\pi_{V_i}f\|^2\\
	&\le 2(\sum_{i\in I}w_i^2\|(a_i\pi_{W_i}-b_i\pi_{V_i})f\|^2+\sum_{i\in I}w_i^2\|b_i\pi_{V_i}f\|^2)\\
	&\le  2(\lambda_1\sum_{i\in I}a_i^2w_i^2\|\pi_{W_i}f\|^2+\lambda_2\sum_{i\in I}b_i^2w_i^2\|\pi_{V_i}f\|^2+\mu\|K^*f\|^2+\sum_{i\in I}w_i^2\|\pi_{V_i}f\|^2),
	\end{align*}
	we have
	$$(1-2\lambda_1)(\inf_{i\in I}a_i)^2\sum_{i\in I}w_i^2\|\pi_{W_i}f\|^2\le 2(1+\lambda_2+\frac{\mu}{A})(\sup_{i\in I}b_i)^2\sum_{i\in I}w_i^2\|\pi_{V_i}f\|^2.$$
	Thus,
	$$\sum_{i\in I}w_i^2\|\pi_{V_i}f\|^2\ge \frac{(1-2\lambda_1)(\inf_{i\in I}a_i)^2}{2(1+\lambda_2+\frac{\mu}{A})(\sup_{i\in I}b_i)^2}\sum_{i\in I}w_i^2\|\pi_{W_i}f\|^2\ge\frac{(1-2\lambda_1)(\inf_{i\in I}a_i)^2}{2(1+\lambda_2+\frac{\mu}{A})(\sup_{i\in I}b_i)^2}A\|K^*f\|^2.$$ 
	Hence $\{(W_i,w_i)\}_{i\in I}$ is a $K$-fusion frame for $\hs$ with bounds
$$\frac{(1-2\lambda_1)(\inf_{i\in I}a_i)^2}{2(1+\lambda_2+\frac{\mu}{A})(\sup_{i\in I}b_i)^2}A~ {\rm and}~\frac{2(1+\lambda_1+\frac{\mu}{A})(\sup_{i\in I}a_i)^2}{(1-2\lambda_2)(\inf_{i\in I}b_i)^2}B.$$
\end{proof}
Next, we study the other result of stability of $K$-fusion frames under perturbation.
\begin{theorem}
	Let $\ws=\{(W_i,w_i)\}_{i\in I}$ be a $K$-fusion frame for $\hs$ with bounds $A$ and $B$. Suppose that $\{V_i\}$ is a family of closed subspaces of $\hs$ and there exist constants $\lambda_1,~\lambda_2,~\mu\ge 0$ such that $\max\{\lambda_1+\frac{\mu}{\sqrt{A}},\lambda_2\}<1$ and
	$$\bigg(\sum_{i\in I}w_i^2\|(\pi_{W_i}-\pi_{V_i})f\|^2\bigg)^{1/2}\le \lambda_1\bigg(\sum_{i\in I}w_i^2\|\pi_{W_i}f\|^2\bigg)^{1/2}+\lambda_2\bigg(\sum_{i\in I}w_i^2\|\pi_{V_i}f\|^2\bigg)^{1/2}+\mu\|K^*f\|^2$$
	for all $f\in\hs$. Then $\{(V_i,w_i)\}_{i\in I}$ is a $K$-fusion frame for $\hs$ with bounds
	$$A\bigg(1-\frac{\lambda_1+\lambda_2+\frac{\mu}{\sqrt{A}}}{1+\lambda_2}\bigg)^2~~~~~{\rm and}~~~~~B\bigg(1+\frac{\lambda_1+\lambda_2+\frac{\mu}{\sqrt{A}}}{1-\lambda_2}\bigg)^2.$$
\end{theorem}
\textbf{Proof.}  Since $$A\|K^*f\|^2\le\sum_{i\in I}w_i^2\|\pi_{W_i}f\|^2,$$ then
$$-\|K^*f\|\ge -\frac{1}{\sqrt{A}}\bigg(\sum_{i\in I}w_i^2\|\pi_{W_i}f\|^2\bigg)^{1/2}.$$
For all $f\in\hs$ we have
\begin{eqnarray}
\bigg(\sum_{i\in I}w_i^2\|\pi_{V_i}f\|^2\bigg)^{1/2}
&=& \bigg(\sum_{i\in I}w_i^2(\|\pi_{V_i}f\|^2-\|\pi_{W_i}f\|^2+\|\pi_{W_i}f\|^2)\bigg)^{1/2}  \nonumber\\
&\ge& \bigg(\sum_{i\in I}w_i^2\|\pi_{W_i}f\|^2\bigg)^{1/2}-\bigg(\sum_{i\in I}w_i^2\|(\pi_{W_i}-\pi_{V_i})f\|^2\bigg)^{1/2}\nonumber\\
&\ge&(1-\lambda_1-\frac{\mu}{\sqrt{A}}) \bigg(\sum_{i\in I}w_i^2\|\pi_{W_i}f\|^2\bigg)^{1/2}-\lambda_2\bigg(\sum_{i\in I}w_i^2\|\pi_{V_i}f\|^2\bigg)^{1/2}.    \nonumber
\end{eqnarray}
Hence
\begin{eqnarray}
\bigg(\sum_{i\in I}w_i^2\|\pi_{V_i}f\|^2\bigg)^{1/2}
&\ge&\bigg(1-\frac{\lambda_1+\lambda_2+\frac{\mu}{\sqrt{A}}}{1+\lambda_2}\bigg)\bigg(\sum_{i\in I}w_i^2\|\pi_{W_i}f\|^2\bigg)^{1/2}  \nonumber\\
&\ge&\sqrt{A}\bigg(1-\frac{\lambda_1+\lambda_2+\frac{\mu}{\sqrt{A}}}{1+\lambda_2}\bigg)\|K^*f\|.    \nonumber
\end{eqnarray}
On the other hand, for all $f\in\hs$ we have
\begin{eqnarray}
\bigg(\sum_{i\in I}w_i^2\|\pi_{V_i}f\|^2\bigg)^{1/2}
&\le& \bigg(\sum_{i\in I}w_i^2\|\pi_{W_i}f\|^2\bigg)^{1/2}+\bigg(\sum_{i\in I}w_i^2\|(\pi_{W_i}-\pi_{V_i})f\|^2\bigg)^{1/2}\nonumber\\
&\le&(1+\lambda_1+\frac{\mu}{\sqrt{A}}) \bigg(\sum_{i\in I}w_i^2\|\pi_{W_i}f\|^2\bigg)^{1/2}+\lambda_2\bigg(\sum_{i\in I}w_i^2\|\pi_{V_i}f\|^2\bigg)^{1/2}.    \nonumber
\end{eqnarray}
Hence
\begin{eqnarray}
\bigg(\sum_{i\in I}w_i^2\|\pi_{V_i}f\|^2\bigg)^{1/2}
&\le&\bigg(1+\frac{\lambda_1+\lambda_2+\frac{\mu}{\sqrt{A}}}{1-\lambda_2}\bigg)\bigg(\sum_{i\in I}w_i^2\|\pi_{W_i}f\|^2\bigg)^{1/2}  \nonumber\\
&\le&\sqrt{B}\bigg(1-\frac{\lambda_1+\lambda_2+\frac{\mu}{\sqrt{A}}}{1+\lambda_2}\bigg)\|f\|.    \nonumber
\end{eqnarray}
The proof is completed.\qed
\section{Atomic systems}
In this section, we study atomic systems for a subspace sequence.
\begin{definition}\label{def5}
	Let $K\in\ls(\hs)$. Let $\ws=\{(W_i,w_i)\}_{i\in I}$ be a sequence of closed subspaces in $\hs$ and let $\{w_i\}_{i\in I}$ be a family of positive weights. Then $\ws$ is called an atomic system for $K$, if the following conditions are satisfied:
	\begin{enumerate}
		\item[(i)] $\ws$ is a Bessel fusion sequence;
		\item [(ii)]For every $f\in\hs$, there exist $\{a_i\}_{i\in I}\in \bigg(\sum_{i\in I}\oplus W_i\bigg)_{\ell^2}$ such that $\|\{a_i\}\|_{(\sum_{i\in I}\oplus W_i)_{\ell^2}}\le C\|f\|$ for some $C>0$ and
		$$Kf=\sum_{i\in I}w_ia_i.$$
	\end{enumerate}
\end{definition}
The following theorem  gives  the existence of the atomic systems for an operator.
\begin{theorem}\cite{guavructa2012frames}
	Let $\hs$ be a separable Hilbert space and $K\in\ls(\hs)$. Then $K$ has an atomic system.
\end{theorem}
Every operator $K$ has an atomic system with a converge sequence. By Remark \ref{rem2} and Definition \ref{def5}, we known that the sequence must be a Bessel sequence. Conversely, one may ask  whether every Bessel fusion sequence $\ws=\{(W_i,w_i)\}_{i\in I}$ has an operator $K$ which makes $\ws$ an atomic system for $K$. The answer is affirmative by the following result.
\begin{theorem}\label{thm12}
	Let $\ws=\{(W_i,w_i)\}_{i\in I}$ be a Bessel fusion sequence in $\hs$. Then $\ws$ is an atomic system for the frame operator $S_{\ws}$.
\end{theorem}
\begin{proof}

 Suppose  $\ws=\{(W_i,w_i)\}_{i\in I}$ is Bessel fusion sequence in $\hs$ with Bessel bound $B$. Let $S_{\ws}$ be the frame operator of $\ws$ and suppose $\{f_{ij}\}_{j\in J_i}$ is a Parseval frame for $W_i$, then we have $\{w_if_{ij}\}_{j\in J_i,i\in I}$ is a frame for $\hs$ with frame operator $S_{\ws}$, and then
$$S_{\ws}f=\sum_{i\in I}w_i^2\pi_{W_i}f=\sum_{j\in J_i}\sum_{i\in I}\left\langle f,w_if_{ij}\right\rangle w_if_{ij},~~~\forall f\in\hs.$$
Hence $S_{\ws}$ is bounded on $\hs$.\\
Let $\{a_i\}_{i\in I}=\{w_i\pi_{W_i}f\}_{i\in I}\in (\sum_{i\in I}\oplus W_i)_{\ell^2}$. We have
$$\|\{a_i\}\|_{(\sum_{i\in I}\oplus W_i)_{\ell^2}}^2=\|\{w_i\pi_{W_i}f\}_{i\in I}\|_{(\sum_{i\in I}\oplus W_i)_{\ell^2}}^2=\sum_{i\in I}w_i^2\|\pi_{W_i}f\|^2\le B\|f\|^2,$$
then $\|\{a_i\}\|_{(\sum_{i\in I}\oplus W_i)_{\ell^2}}^2\le B\|f\|^2$. Hence, $\|\{a_i\}\|_{(\sum_{i\in I}\oplus W_i)_{\ell^2}}\le \sqrt{B}\|f\|$ for each $f\in\hs$. From the Definition 5, $\ws$ is an atomic system for the frame operator $S_{\ws}$. 
\end{proof}

 In fact, by Theorem \ref{thm12}, we can find that every Bessel sequence has an operator $K$ which makes it an atomic system for $K$.
Next, we give a characterization of atomic systems with subspace sequences.
\begin{theorem}\label{thm13}
	Let $\{(W_i)\}_{i\in I}$ be a family of closed subspaces of $\hs$ and $\{w_i\}_{i\in I}$ be a family of positive weights.  Then the following statements are equivalent
	\begin{enumerate}
		\item[(i)] $\ws=\{(W_i,w_i)\}_{i\in I}$ is an atomic system for $K$;
		\item[(ii)] $\ws=\{(W_i,w_i)\}_{i\in I}$ is a $K$-fusion frame for $\hs$, i.e., there exist constant $0<A\le B<\infty$ such that
		$$	A\|K^*f\|^2\le\sum_{i\in I}w_i^2\|\pi_{W_i}(f)\|^2\le B\|f\|^2,~~~\forall f\in\hs.$$
	\end{enumerate}
\end{theorem}
\begin{proof} (i)$\Rightarrow$(ii) Since $\ws=\{(W_i,w_i)\}_{i\in I}$ is an atomic system for $K$, $\ws$ is a Bessel fusion frame for $\hs$. We now prove that there exists a constant $A>0$ such that
$$\sum_{i\in I}w_i^2\|\pi_{W_i}f\|^2\ge A\|K^*f\|^2$$
for all $f\in\hs$.
By taking $g$ instead of $f$ in the condition (ii) of Definition \ref{def5}, there exists $c>0$ such that for every $g\in\hs$ there exists $\{a_i\}_{i\in I}\in \bigg(\sum_{i\in I}\oplus W_i\bigg)_{\ell^2}$ such that $\|a_f\|_{(\sum_{i\in I}\oplus W_i)_{\ell^2}}\le C\|f\|$ and
$$Kg=\sum_{i\in I}w_ia_i=\sum_{i\in I}w_i\pi_{W_i}a_i.$$
By using Cauchy-Schwarz inequality for all $f\in\hs$ we have
\begin{eqnarray}
\|K^*f\|^2&=&\sup_{\|g\|=1}|\left\langle K^*f,g\right\rangle |^2=\sup_{\|g\|=1}|\left\langle f,Kg\right\rangle |^2 \nonumber\\
&=& \sup_{\|g\|=1}|\langle f,\sum_{i\in I}w_i\pi_{W_i}a_i\rangle | ^2=\sup_{\|g\|=1}|\sum_{i\in I}\langle w_i\pi_{W_i}f,a_i\rangle | ^2 \nonumber\\
&\le& \sup_{\|g\|=1}\big( \sum_{i\in I}|a_i|^2\big) \cdot\big( \sum_{i\in I}w_i^2\|\pi_{W_i}f\|^2\big)\nonumber\\
&\le&c \sum_{i\in I}w_i^2\|\pi_{W_i}f\|^2.   \nonumber
\end{eqnarray}
By setting $A=\frac{1}{C}$, we get the desired result.

(ii)$\Rightarrow$(i) Let $A,B$ be the bounds of $\ws$. We know that $\ws$ is a fusion Bessel sequence. Hence the synthesis operator $T_{\ws}$ is bounded and we have
$$A\|K^*f\|^2\le \sum_{i\in I}w_i^2\|\pi_{W_i}f\|^2=\|T^*_{\ws}f\|^2,$$
for all $f\in\hs$. This implies that
$$KK^*\le \frac{1}{A}T_{\ws}T_{\ws}^*.$$
Then Theorem \ref{thm3} implies that there exists a bounded operator $\Gamma:\hs\rightarrow (\sum_{i\in I}\oplus W_i)_{\ell^2}$ such that $K=T_{\ws}\Gamma$. For every $f\in\hs$, we define $\Gamma f=\{a_i\}_{i\in I}$. Then we have
$$Kf=T_{\ws}\Gamma f=T_{\ws}\{a_i\}_{i\in I}=\sum_{i\in I}w_ia_i$$
and
$$\|\{a_i\}\|_{(\sum_{i\in I}\oplus W_i)_{\ell^2}}=\|\Gamma f\|_{(\sum_{i\in I}\oplus W_i)_{\ell^2}}\le \|\Gamma\|\|f\|$$
for all $f\in\hs$. From the Definition \ref{def5}, it shows that $\ws$ is an atomic system for $K$.
\end{proof}
\begin{corollary}
	Let $K\in\ls(\hs)$ and let $\{e_i\}_{i\in I}$ be a orthonormal basis for $\hs$, then $\{(Ke_i,1)\}_{i\in I}$ is a $K$-fusion frame for $\hs$.
\end{corollary}
\begin{proof} By Proposition 4.6 of \cite{balazs2011classification},  the Bessel sequence for $K$ are precisely the families $\{Ke_i\}_{i\in I}$. By Theorem \ref{thm12}, $\{(Ke_i,1)\}_{i\in I}$ is a atomic system. Then by Theorem \ref{thm13}, the conclusion is true.
\end{proof}

Let $U$ and $V$ be bounded (linear) operators on a Hilbert  space $\hs$ with the kernel condition
$N(V)\subset N(U)$,
then the quotient $[U/V]$ is a map from $R(V)$ to $R(U)$
defined by $Vf\longmapsto Uf$ for all $f\in\hs$. We note that $P=[U/V] $ is a linear operator on $\hs$ if and only if $N(V)\subset N(U)$. In this case $D(P)=R(V)$, $R(P)\subset R(U)$ and $PV=U$.  
The quotient $[U/V]$ is called a semiclosed operator and its collection is closed under sum and product \cite{Kaufman1979Semiclosed}. 

In the following theorem we generalize Theorem 6 of \cite{li2017some}(or Theorem 5.1 of \cite{ramu2016frame}) to $K$-fusion frames.
\begin{theorem}
	Let $K\in\hs$ and $\ws=\{(W_i,w_i)\}_{i\in I}$ be a Bessel  fusion sequence in $\hs$ with the frame operator $S_{\ws}$.Then $\ws=\{(W_i,w_i)\}_{i\in I}$ is an  atomic system for $K$ if and only if the quotient operator $[K^*/S^{1/2}_{\Lambda}]$ is bounded.
\end{theorem}

\begin{proof}
	$\Longrightarrow:$ Since $\ws=\{(W_i,w_i)\}_{i\in I}$ is an  atomic system for $K$, from Theorem \ref{thm13}, there exists a constant $A>0$ such that 
	$$A\|K^*f\|^2\le \sum_{i\in I}w_i^2\|\pi_{W_i}f\|^2=\left\langle S_{\ws}f,f\right\rangle ,~~~\forall f\in\hs.$$
	That is, $A\|K^*f\|^2\le \|S_{\ws}^{1/2}f\|^2$ for all $f\in\hs$.
	Define $P:R(S_{\ws}^{1/2})\rightarrow R(K^*)$ by
	$$P(S_{\ws}^{1/2}f)=K^*f,~~\forall f\in\hs.$$
	Then $P$ is well-defined because $N(S_{\ws}^{1/2})\subset N(K^*)$. For all $f\in\hs$, we have
	$$\|PS_{\ws}^{1/2}f\|=\|K^*f\|\le \frac{1}{\sqrt{A}}\|S_{\ws}^{1/2}f\|.$$
	So $P$ is bounded. From the notion of quotient of bounded operators, $P$ can be expressed as  $[K^*/S^{1/2}_{\ws}]$.
	
	$\Longleftarrow:$ Suppose that the quotient operator   $[K^*/S^{1/2}_{\ws}]$ is bounded. Then there exists $\lambda>0$ such that
	$$\|K^*f\|^2\le \lambda\|S_{\ws}^{1/2}f\|^2,~~~\forall f\in\hs.$$
	Thus
	$$\frac{1}{\lambda}\|K^*f\|^2\le\|S_{\ws}^{1/2}f\|^2=\left\langle S_{\ws}f,f\right\rangle =\sum_{i\in I}w_i^2\|\pi_{W_i}f\|^2,$$
	for all $f\in\hs$. Hence   $\{(W_i,w_i)\}_{i\in I}$ is a $K$-fusion frame for $\hs$. From Theorem \ref{thm13}, the result holds.
\end{proof}

In the following proposition by extending Theorem 6 of \cite{li2017some}, we characterize the equivalence of atomic systems and quotient operators.
\begin{theorem}
	Let $\ws=\{(W_i,w_i)\}_{i\in I}$ be an atomic system for $K$ with frame operator $S_{\ws}$. Let $T\in\ls(\hs)$ has closed range with $T^{\dagger}T(W_i)\subset W_i$.  Then the following are equivalent:
	\begin{enumerate}
		\item [(i)]$\{(TW_i,w_i)\}_{i\in I}$ is an atomic system for $TK$;
		\item [(ii)]$[(TK)^*/S^{1/2}_{\ws}T^*]$ is bounded;
		\item [(iii)]$[(TK)^*/(TS_{\ws}T^*)^{1/2}]$ is bounded.
	\end{enumerate}	
\end{theorem}
\begin{proof}
	(i)$\Rightarrow$(ii); From Theorem \ref{thm13},  $\{\Lambda_iT\}_{i\in I}$ is a $TK$-fusion frame. Then there exist $\lambda>0$ such that
	\begin{equation}\label{3.01}
	\lambda\|(TK)^*f\|^2\le \sum_{i\in I}w_i^2\|\pi_{TW_i}f\|^2
	\end{equation}
	from \eqref{2.7}, we have
		\begin{equation}\label{3.02}
\sum_{i\in I}w_i^2\|\pi_{TW_i}f\|^2\le\|{T^{\dagger}}^{*}\|^2 \sum_{i\in I}w_i^2\|\pi_{W_i}T^*f\|^2=\|{T^{\dagger}}^{*}\|^2\left\langle S_{\ws}T^*f,T^*f\right\rangle =\|{T^{\dagger}}^{*}\|^2\|S^{1/2}_{\ws}T^*f\|^2,~~\forall f\in\hs.
	\end{equation}
By \eqref{3.01} and \eqref{3.02}, we have
$$\lambda\|{T^{\dagger}}^{*}\|^{-2}\|(TK)^*f\|^2\le \|S^{1/2}_{\ws}T^*f\|^2.$$
	Hence $[(TK)^*/S^{1/2}_{\ws}T^*]$ is bounded.
	
	(ii)$\Rightarrow$(iii); Suppose $[(TK)^*/S^{1/2}_{\ws}T^*]$ is bounded. Then there exists $\mu>0$ such that 
	$$\|(TK)^*f\|^2\le\mu\|S^{1/2}_{\ws}T^*f\|^2,~~\forall f\in\hs.$$
	Since
	\begin{align}
	\|(TS_{\ws}T^*)^{1/2}f\|^2&=\left\langle (TS_{\ws}T^*)^{1/2}f,(TS_{\ws}T^*)^{1/2}f\right\rangle =\left\langle (TS_{\ws}T^*)f,f\right\rangle\nonumber\\
	&=\left\langle  S_{\ws}T^*f,T^*f\right\rangle =\|S^{1/2}_{\ws}T^*f\|^2,
	\end{align}
	for all $f\in\hs$, we have 
	$$\frac{1}{\mu}\|(TK)^*f\|^2\le \|(TS_{\ws}T^*)^{1/2}f\|^2.$$
	Therefore $[(TK)^*/(TS_{\ws}T^*)^{1/2}]$ is bounded.
	
	(iii)$\Rightarrow$(i); Suppose $[(TK)^*/(TS_{\ws}T^*)^{1/2}]$ is bounded, then there exists $\mu>0$ such that 
	$$\|(TK)^*f\|^2\le \mu\|(TS_{\ws}T^*)^{1/2}f\|^2,~~~\forall f\in\hs.$$
	From \eqref{2.5}, we have
	$$\sum_{i\in I}w_i^2\|\pi_{TW_i}f\|^2\ge \|T^*\|^{-2}\sum_{i\in I}w_i^2\|\pi_{W_i}T^*f\|^2=\|T^*\|^{-2}\left\langle  S_{\ws}T^*f,T^*f\right\rangle=\|T^*\|^{-2}\left\langle  TS_{\ws}T^*f,Tf\right\rangle.$$
	So $TS_{\ws}T^*$ is positive and self-adjoint, its square root exists, and it is denoted by $(TS_{\ws}T^*)^{1/2}$. Hence 
	$$\sum_{i\in I}w_i^2\|\pi_{TW_i}f\|^2\ge\|(TS_{\ws}T^*)^{1/2}f\|^2\ge\frac{1}{\mu}\|(TK)^*f\|^2,~~\forall f\in\hs.$$
	
	On the other hand, from \eqref{2.7},
	$$\sum_{i\in I}w_i^2\|\pi_{TW_i}f\|^2\le \|{T^{\dagger}}^{*}\|^2 \sum_{i\in I}w_i^2\|\pi_{W_i}T^*f\|^2\le B\|{T^{\dagger}}^{*}\|^2\|T^*\|^2\|f\|^2,$$
	where $B$ is the upper frame bound of $\ws$. 
	Hence 	 $\{(TW_i,w_i)\}_{i\in I}$ is a $TK$-fusion frame, and from Theorem \ref{thm13}, as desired.
\end{proof}
In the following result, we give a method to construct $K$-fusion frame in $\hs$ from the view of atomic systems.
\begin{theorem}
	Let $K\in\ls(\hs)$ be with closed range and let $K^{\dagger}$ be the pseudo-inverse of $K$ and satisfy $K^{\dagger}K(W_i)\subset W_i$. If $\{(W_i,w_i)\}_{i\in I}$ is a fusion frame for $\hs$, then $\ws=\{(KW_i,w_i)\}_{i\in I}$ is a $K$-fusion frame for $\hs$.
\end{theorem}
\begin{proof}
 By Theorem \ref{thm13}, we just need to show that $\{(W_i,w_i)\}_{i\in I}$ is an atomic system for $K$. \\
Suppose $\{(W_i,w_i)\}_{i\in I}$ is a fusion frame for $\hs$ with frame operator $S_{\ws}$, for all $f\in\hs$, we have
\begin{equation}\label{3.1}
S_{\ws}f=\sum_{i\in I}w^2_i\pi_{W_i}f.
\end{equation}
With $Kf$ instead of $f$ in \eqref{3.1}, we have
$$S_{\ws}Kf=\sum_{i\in I}w^2_i\pi_{W_i}Kf.$$
Hence,
$$Kf=\sum_{i\in I}w_i^2S_{\ws}^{-1}\pi_{W_i}Kf,~~\forall f\in\hs.$$
Then we show that $\{(KW_i,w_i)\}_{i\in I}$ is a Bessel fusion sequence and $\|\{w_iS_{\ws}^{-1}\pi_{W_i}f\}\|_{(\sum_{i\in I}\oplus W_i)_{\ell^2}}\le C\|f\|$, where $C$ is a positive constant. \\
Let $A,B$ be the lower and upper frame bounds of $\ws$, then for any $f\in\hs$, we have
$$\sum_{i\in I}w_i^2\|S_{\ws}^{-1}\pi_{W_i}f\|^2\le A^{-1}\|f\|^2.$$
On the other hand, from Lemma \ref{lem2}, we obtain, with $K^{\dagger}$ instead of $K$:
$$\|\pi_{KW_i}f\|\le\|{K^{\dagger}}^{*}\|\|\pi_{W_i}K^*f\|.$$
Then
\begin{eqnarray}
\sum_{i\in I}w_i^2\|\pi_{KW_i}f\|^2
&\le& \|{K^*}^{\dagger}\|^2\sum_{i\in I}w_i^2\|\pi_{W_i}K^*f\|^2  \nonumber\\
&\le& \|{K^*}^{\dagger}\|^2B\|K^*f\|^2\nonumber\\
&\le&B \|{K^*}^{\dagger}\|^{2}\|K\|^2 \|f\|^2.   \nonumber
\end{eqnarray}
So $\{(W_i,w_i)\}_{i\in I}$ is an atomic system for $K$. By Theorem \ref{thm13} the conclusion holds.
\end{proof}
\section*{Acknowledgements}
The research is supported by the National Natural Science Foundation of China (LJT10110010115). 
\bibliographystyle{plain}
\bibliography{kframe}
\end{document}